\documentclass[oneside]{amsart}

\theoremstyle{plain}    
\newtheorem{thm}{Theorem}[section]
\numberwithin{equation}{section} 
\numberwithin{figure}{section} 
\theoremstyle{plain}    
\newtheorem*{thm*}{Theorem} 
\theoremstyle{plain}    
\newtheorem{cor}[thm]{Corollary} 
\theoremstyle{plain}    
\newtheorem{prop}[thm]{Proposition} 
\theoremstyle{definition}
\newtheorem*{defn*}{Definition}
\theoremstyle{remark}
\newtheorem{rem}[thm]{Remark}
\newenvironment{lyxlist}[1]
  {\begin{list}{}
    {\settowidth{\labelwidth}{#1}
     \setlength{\leftmargin}{\labelwidth}
     \addtolength{\leftmargin}{\labelsep}
     }}
  {\end{list}}

\begin{document}

\title[Quantum periods--I]{Quantum periods -- I. Semi-infinite variations of Hodge structures.} \thanks{Preprint ENS DMA-00-19, June 7, 2000}

\author{Serguei Barannikov}

\maketitle
\centerline{\'{E}cole Normale Sup\'{e}rieure, 45 rue d'Ulm, 75230 Paris, CNRS
UMR8553} \centerline{e-mail: serguei.barannikov@ens.fr}

\urladdr{http://www.dma.ens.fr/\~{}barannik/}

\begin{abstract}
We introduce a generalization of variations of Hodge structures living over
moduli spaces of non-commutative deformations of complex manifolds. Hodge structure
associated with a point of such moduli space is an element of Sato type grassmanian
of semi-infinite subspaces in \break \( H^{*}(X,\mathbb {C})[[\hbar ^{-1},\hbar ]] \).
Periods associated with such semi-infinite Hodge structures serve in order to
extend mirror symmetry relations in dimensions greater then three. 
\end{abstract}
\setcounter{tocdepth}{1}

\tableofcontents{}

\section{Introduction\label{intr}}

This paper is concerned with theory of periods of non-commutative varieties
which arise naturally in the context of mirror symmetry.

We have noticed in \cite{B2} that in order to understand mirror symmetry phenomena
in dimensions higher then three one has to enlarge the class of functions provided
by periods of complex manifolds. In this paper we introduce a new class which
consists of periods associated with \emph{non-commutative} deformations of complex
manifolds. We show that this is precisely the class which is needed in order
to generalize to dimensions higher then three the mirror symmetry relation between
rational Gromov-Witten invariants of Calabi-Yau manifolds and periods associated
with their mirror partners.

As another motivation for developing theory of quantum periods let us mention
that such a theory should serve as one of guides for future theory of non-commutative
varieties. It is interesting to make a somewhat ambitious parallel and to note
that, historically, developement of integration theory of analytic functions
(theory of elliptic integrals, in particular, and, more generally, of Abelian
integrals) was one of the principal reasons for introduction of Riemann surfaces.

One may hope also that quantum periods should help to understand certain aspects
of varieties from ordinary commutative world (see, for example, the local Torelli
property of the period map from \S \ref{permapM}) and, in particular, especially
those related with \( K- \)theory.

We plan to write a series of papers concerning the theory of quantum periods
and its applications. Next paper (\cite{B3}) is devoted to applications to
non-commutative geometry and homological mirror symmetry. In a separate paper
(\cite{B4}) we consider in details an oversimplified situation where periods
of somewhat related type appear. There we deal with semi-infinite Hodge structure
associated with mirror partners of projective spaces. The main result of \cite{B4}
is the equality expressing generating functions for total sets of rational Gromov-Witten
invariants of \( \mathbb {C}\mathbb {P}^{n} \)in terms of periods of semi-infinite
Hodge structures associated with their mirrors.

Let us say a few words on the plan of the paper. Let \( X \) denotes a complex
smooth manifold. We begin in \S \ref{section m.spaces} by recalling description
of the moduli space \( \mathcal{M} \) of non-commutative deformations of \( X \)
via differential graded Lie algebra \( \mathfrak {g}_{A_{\infty }} \) of infinitesimal
\( A_{\infty }- \)symmetries of \( X \). The moduli space \( \mathcal{M} \)
is the quotient of graded scheme of zeroes of solutions to Maurer-Cartan equation
in \( \mathfrak {g}_{A_{\infty }} \) with respect to gauge equivalences. In
\S \ref{gm} we show that Gauss-Manin local system extends naturally over the
moduli space \( \mathcal{M} \). This is based on the general observation (prop.\ref{pr:dmod})
that sheaves with flat connection over moduli space associated with differential
graded Lie algebra \( \mathfrak {g} \) can be described via modules over \( \mathfrak {g}\otimes \mathbb {C}[\xi ]/\xi ^{2} \).
In \S \ref{sect:semiHodge} we associate with moduli space \( \mathcal{M} \)
semi-infinite analog of variations of Hodge structures. If \( \gamma \in (\mathfrak {g}_{A_{\infty }}\otimes \mathfrak {M}_{R})^{1} \)
denotes a solution to Maurer-Cartan equation describing a non-commutative deformation
of \( X \) given by \( \phi \in \mathcal{M}(R) \) then the semi-infinite analog
of Hodge filtration associated with \( \phi  \) is the subspace \( L(\phi )\in Gr_{\frac{\infty }{2}}(H[[\hbar ^{-1},\hbar ]])(R) \)
which consists of elements of the form 
\begin{equation}
\left[ l_{\hbar }\exp (\frac{1}{\hbar }i_{\gamma })(\varphi _{0}+\varphi _{1}\hbar +\varphi _{2}\hbar ^{2}\ldots )\right] \: \: \text {where}\, \varphi _{i}\in \Omega ^{*}(X,\mathbb {C})\otimes R\notag 
\end{equation}
 where \( i_{\gamma } \) denotes contraction \( \Lambda ^{*}T\otimes \Omega ^{*}\rightarrow \Omega ^{*} \)extended
by wedge multiplication by \( (0,*)- \)forms and \( l_{\hbar } \) denotes
rescaling operator \( \varphi ^{p,q}\rightarrow \hbar ^{\frac{q-p}{2}}\varphi ^{p,q} \)
for \( \varphi ^{p,q}\in \Omega ^{p,q} \). In \S \ref{section griffiths} we
study variations of the subspaces \( L(\phi ) \), \( \phi \in \mathcal{M}(R) \)
and prove that they satisfy semi-infinite version of Griffiths transversality
condition: 
\[
\nabla L\subset \hbar ^{-1}L\]
 In \S \ref{mirror} we explain the role played by periods of such semi-infinite
Hodge structures in higher dimensional mirror symmetry. Similar to the three-dimensional
case (\cite{CGOP}) the calculation of Gromov-Witten invariants of an \( n- \)dimensional
Calabi-Yau manifold \( Y \) via quantum periods associated with mirror Calabi-Yau
manifold \( X \) is done in three steps. At the first step one specifies normalization
of the periods. Assume that we are given a linear subspace \( S\subset H[[\hbar ^{-1},\hbar ]] \),
\( \hbar ^{-1}S\subset S \) transversal to subspace \( L(0) \) and an element
\( \psi _{0}\in L(0) \). Such linear subspace \( S=L_{W} \) and element \( \psi _{0} \)
are associated naturally with monodromy weight filtration \( W \) of maximal
unipotency cusp on moduli space of complex structures on \( X \). Then there
exists unique element \( \Psi ^{W}(t,\hbar )\in L(t) \), \( t\in \mathcal{M} \)
satisfying the normalization condition \( \Psi ^{W}(t,\hbar )\in (L_{W}+\psi _{0}) \).
Second step is a choice of coordinates on the moduli space \( \mathcal{M} \)
using the periods normalized in the first step. The coordinates are induced
via map to linear space \( \mathcal{M}\rightarrow L_{W}/\hbar ^{-1}L_{W} \)
defined as \( (\Psi ^{W}(t,\hbar )-\psi _{0})\, \text {mod}\, \hbar ^{-1}L_{W} \).
In the last step one writes differential equation: 
\[
\frac{\partial ^{2}\Psi ^{W,\alpha }}{\partial t^{a}_{W}\partial t^{b}_{W}}=\hbar ^{-1}\sum _{c}A_{ab}^{c}(t_{W})\frac{\partial \Psi ^{W,\alpha }}{\partial t_{W}^{c}}\]
 satisfied by quantum periods \( \Psi ^{W,\alpha }(t_{W})=\int _{\Delta _{\alpha }}\Psi ^{W} \)
where \( \{\Delta _{\alpha }\} \) is a basis in \( \oplus _{k}H_{k}(X,\mathbb {Z}) \)
and \( \{t_{W}^{a}\} \) are the coordinates from step two. Then \( A_{ab}^{c}(t_{W}) \)
should coincide with structure constants of (big) quantum cohomology of \( Y \)
and the coefficients of appropriate Taylor expansion of \( A_{ab}^{c}(t_{W}) \)
should give the set of all rational Gromov-Witten invariants of \( Y \). We
checked this in \cite{B2} \S 6 for Calabi-Yau complete intersections varieties.

\section{Moduli spaces of \protect\( A_{\infty }-\protect \)deformations of complex
manifolds.\label{section m.spaces}}

It is explained here how to think about moduli space \( \mathcal{M} \) of non-commutative
deformations of complex manifold \( X \). We assume that the reader is familiar
with the technique of deformation theory (see for example, \cite{K1} \S 4,
\cite{B2} \S 2, and references therein). Idea of considering \( A_{\infty }- \)deformations
of \( D^{b}Coh(X) \) was suggested in \cite{K2}. The definition of the moduli
space given below appeared in \cite{BK}.

We would like to stress that a satisfactory definition of objects that are parameterized
by moduli space \( \mathcal{M} \) remains yet a mystery. The only thing which
is known is that it should be some quantum varieties such that their equivalence
classes are in one-to-one correspondence with equivalence classes of \( A_{\infty }- \)deformations
of bounded derived category of coherent sheaves on \( X \). Remark that, as
it follows easily from the definition of \( \mathcal{M} \) given below, our
moduli space contains a closed subspace which parameterizes classes of equivalences
of Poisson brackets on \( X \). The space of equivalence classes of first-order
deformations parameterized by \( \mathcal{M} \) is given by \( T_{[X]}\mathcal{M}=\oplus _{p,q}H^{q}(X,\Lambda ^{p}T_{X})[p-q] \).
The subspace of Poisson brackets corresponds on the level of first-order deformations
to the subspace \( H^{0}(X,\Lambda ^{2}T_{X})\subset (\oplus _{p,q}H^{q}(X,\Lambda ^{p}T_{X})[p-q]) \).

\subsection{Deformations of complex structure}

Recall first the standard description of deformations of complex structures.
Given complex manifold \( X \) let \( J \) denotes the corresponding complex
structure on the underlying \( C^{\infty }- \)manifold \( X_{C^{\infty }} \).
According to Kodaira-Spencer theory the deformations of \( J \) are described
by elements \( \rho \in \Omega ^{0,1}(X,T) \) satisfying 
\begin{equation}
\label{mc}
\overline{\partial }\rho +\frac{1}{2}[\rho ,\rho ]=0
\end{equation}
 Namely, to an element \( \rho  \) one can associate a deformation of the decomposition
\( T_{\mathbb {C}}^{*}=T^{*}\bigoplus \overline{T}^{*} \)of bundle of \( \mathbb {C}- \)valued
\( 1- \)forms into the sum of sub-bundles corresponding to differentials of
holomorphic and anti-holomorphic coordinates. The deformed sub-bundle \( T^{*} \)
is given by the graph of the linear map corresponding to \( \rho \in \Gamma (X,T\bigotimes \overline{T}^{*}) \).
Eq. (\ref{mc}) means that this sub-bundle defines formally integrable distribution.
Such a distribution gives a new complex structure on \( X_{C^{\infty }} \)
by Newlander-Nirenberg theorem. Also, equivalent deformations are related via
the action of the Lie algebra \( \Gamma _{C^{\infty }}(X,T) \). Notice that
this description is formulated intrinsically in terms of of differential graded
Lie algebra on \( \oplus _{q}\Omega ^{0,q}(X,T)[-q] \). We denote this differential
graded Lie algebra via \( \frak {g}_{cs} \).

\subsection{\protect\( A_{\infty }-\protect \)deformations}

The \( A_{\infty }- \)deformations of \( X \) allow similar description using
elements \( \gamma \in \Omega ^{0,*}(X,\Lambda ^{*}T) \). Consider the differential
graded Lie algebra 
\[
\frak {g}_{A_{\infty }}:=\bigoplus _{i}\frak {g}_{A_{\infty }}^{i}[-i],\: \frak {g}_{A_{\infty }}^{i}=\bigoplus _{q-p+1=i}\Omega ^{0,q}(X,\Lambda ^{p}T)\]
 equipped with Dolbeault differential and the bracket which is given by Schouten-Nijenhuis
bracket on \( \Lambda ^{*}T \) extended using cup-product of anti-holomorphic
differentials. Consider the graded scheme of solutions to Maurer-Cartan equation
\[
MC_{A_{\infty }}:=\{\overline{\partial }\gamma +\frac{1}{2}[\gamma ,\gamma ]=0\}\subset \frak {g}_{A_{\infty }}[1]\]
 The standard approach consists of developing solutions to Maurer-Cartan in
power series \( \gamma =\sum _{a}\gamma _{a}t^{a}+\frac{1}{2}\sum _{a_{1}a_{2}}\gamma _{a_{1}a_{2}}t^{a_{1}}t^{a_{2}}+\ldots  \).
In this paper we leave aside the question of convergence of power series representing
the solutions, which can be treated separately using the standard technique
of Kuranishi spaces. Consequently, we will work on the level of formal power
series and will be dealing with the completion of \( MC_{A_{\infty }} \) at
zero. In the sequel \( MC_{A_{\infty }} \) will always denote the corresponding
formal scheme. This scheme has points defined over arbitrary graded \( \mathbb {C}- \)Artin
algebra\footnote{%
Recall that any such algebra is isomorphic to \( \mathbb {C}[t_{s}]_{s\in S}/I \),
where \( S \) is some finite set and \( I \) is an ideal containing sufficiently
large power of the maximal ideal \( t\mathbb {C}[t_{s}] \). 
}\( R \): 
\[
MC_{A_{\infty }}(R)=\{\gamma |\overline{\partial }\gamma +\frac{1}{2}[\gamma ,\gamma ]=0,\, \gamma \in (\frak {g}_{A_{\infty }}[1]\otimes \frak {M}_{R})^{0}\}\]
 where \( \frak {M}_{R} \) denotes the unique maximal ideal in \( R \). We
will also consider points of \( MC_{A_{\infty }} \) over projective limits
of Artin algebras\footnote{%
Here and below by Artin algebra we mean graded \( \mathbb {C}- \)Artin algebra
if it does not specified otherwise. 
}, for example, over algebra of formal power series on some finite-dimensional
\( \mathbb {Z}- \)graded vector space. Algebra Lie \( \mathfrak {g}_{A_{\infty }} \)
acts on the scheme \( MC \) via gauge transformations (again it should be understand
as action in the category of formal graded schemes): 
\begin{equation}
\label{gge-trnsf}
\alpha \in \frak {g}_{A_{\infty }}\rightarrow \Dot {\gamma }=\overline{\partial }\alpha +[\gamma ,\alpha ]
\end{equation}
 For any Artin algebra \( R \) the algebra Lie \( (\frak {g}_{A_{\infty }}\otimes \frak {M}_{R})^{0} \)is
nilpotent and the associated group \( \exp \bigl ((\frak {g}_{A_{\infty }}\otimes \frak {M}_{R})^{0}\bigr ) \)
acts on \( MC_{A_{\infty }}(R) \). Explicitly, the action of an element \( \varphi =\exp \alpha \in \exp \bigl ((\frak {g}_{A_{\infty }}\otimes \frak {M}_{R})^{0}\bigr ) \)
is given by 
\begin{equation}
\label{gaugegroup}
\gamma \, \rightarrow \, \gamma ^{\varphi }:=\overline{\partial }\varphi \cdot \varphi ^{-1}+\varphi \cdot \gamma \cdot \varphi ^{-1}=\sum ^{\infty }_{n=0}(\frac{(ad\, \alpha )^{n+1}\overline{\partial }\alpha }{(n+1)!}+\frac{(ad\, \alpha )^{n}\gamma }{n!})
\end{equation}

\begin{defn*}
The (formal) moduli space \( \mathcal{M} \) is the quotient of \( MC_{A_{\infty }} \)by
the gauge transformations (\ref{gaugegroup}). 
\end{defn*}
More precisely this should be understand as follows: \( \mathcal{M} \) is described
via the functor \( {Def}^{\mathbb {Z}}_{\frak {g}_{A_{\infty }}} \): \( \mathbb {Z}- \)graded
Artin algebras \( \rightarrow  \)sets, which is given by 
\begin{equation}
\label{def(R)}
{Def}^{\mathbb {Z}}_{\frak {g}_{A_{\infty }}}(R):=MC_{A_{\infty }}(R)/\exp \bigl ((\frak {g}_{A_{\infty }}\otimes \frak {M}_{R})^{0}\bigr )
\end{equation}
 Morally the value of this functor \( {Def}^{\mathbb {Z}}_{\frak {g}_{A_{\infty }}}(R) \)
on Artin algebra \( R \) is the set of morphisms \( (\text {Spec}\, R,0)\rightarrow \mathcal{M} \).
The tangent space to \( \mathcal{M} \) at the base point (more precisely the
space of equivalence classes of first-order deformations) is equal to \( {Def}^{\mathbb {Z}}_{\frak {g}_{A_{\infty }}}(\mathbb {C}[t]/t^{2}) \)
which gives 
\[
T_{[X]}\mathcal{M}=\oplus _{p,q}H^{q}(X,\Lambda ^{p}T)[p-q]\]
 In the sequel we will denote via \( \mathcal{M}(R) \) the set of \char`\"{}points\char`\"{}
of \( \mathcal{M} \) over \( R \), i.e. set of equivalence classes of \( A_{\infty }- \)deformations
depending on parameters from \( R \).

If one restricts oneself to solutions to Maurer-Cartan equation from \( MC_{cs}(R):=\{\gamma |\overline{\partial }\gamma +\frac{1}{2}[\gamma ,\gamma ]=0,\, \gamma \in \Omega ^{0,1}(X,T)\otimes \frak {M}_{R^{0}}\} \)
then it corresponds to deformations of complex structure on \( X \). We denote
by \( \mathcal{M}^{cs} \) the quotient of \( MC_{cs} \) by gauge equivalences
(\ref{gge-trnsf}). It is easy to see that such equivalences are generated only
by elements from \( \frak {g}_{cs}^{0}=\Gamma _{C^{\infty }}(X,T) \). For compact
complex manifold \( X \) the formal moduli space \( \mathcal{M}_{cs} \) is
isomorphic to the completion at the base point \( [X] \) of Kuranishi space
which is a mini-versal family of deformations of complex structure on \( X \)
(Goldman-Millson) (here we assumed for simplicity that \( H^{0}(X,T_{X})=0 \),
otherwise the statement is slightly more complicated).

The formality theorem from \cite{K1} can be used to show that the moduli space
\( \mathcal{M} \) parameterizes \( A_{\infty }- \)deformations of the category
\( D^{b}Coh(X) \). The corresponding arguments are sketched in \S 3.3 of \cite{B2}.

There is a useful point of view on the quotients of type (\ref{def(R)}) (see
\cite{K1} \S 4, \cite{B2} \S 2). Such quotients can be more generally associated
with germs of \( Q- \)manifolds (or dg-manifolds in different terminology).
This is a germ of \( \mathbb {Z}- \)graded manifold plus a degree one vector
field \( Q \) satisfying \( [Q,Q]=0 \). Any differential graded Lie algebra
can be viewed as \( Q- \)manifold. Namely, vector field 
\[
Q_{\mathfrak {g}}(\gamma )=d\gamma +\frac{1}{2}[\gamma ,\gamma ]\]
 on \( \frak {g}[1] \) satisfies \( [Q_{\mathfrak {g}},Q_{\mathfrak {g}}]=0 \).
Then the moduli space described by \( \mathfrak {g} \) can be viewed as a kind
of nonlinear cohomology associated with \( Q_{\mathfrak {g}} \). More precisely,
it is realized as a quotient 
\begin{equation}
\label{quotent}
\mathcal{M}_{\frak {g}}=(\text {zeros}\, \text {of}\, Q)/\text {distribution}\, \text {generated}\, \text {by}\, \text {vector}\, \text {fields}\, \, \text {of}\, \text {the}\, \text {form}\, [Q,\, \cdot \, ])
\end{equation}
 of the type similar to (\ref{def(R)}). It is often useful to express statements
about the quotient in \( Q- \)equivariant terms on the underlying \( Q \)-manifold.

Let us add that it is natural to expect that for complex algebraic manifold
\( X \) the moduli space \( \mathcal{M} \) can be equipped with appropriate
algebraic structure i.e. \( \mathcal{M} \) must be defined canonically on the
level of (a kind of) dg-scheme.

\section{Gauss-Manin local system over \protect\( \mathcal{M}\protect \)\label{gm}.}

This section is devoted to an analog of Gauss-Manin local system living over
the moduli space of \( A_{\infty }- \)deformations of \( X \).

Sheaves over moduli spaces are described by modules over corresponding differential
graded Lie algebras. The easiest way to see this is to use the geometric language
of \( Q- \)manifolds. Let \( \mathfrak {g} \) denotes differential graded
Lie algebra describing moduli space \( \mathcal{M}_{\mathfrak {g}} \) as in
(\ref{def(R)}). Given dg-module \( h \) over \( \frak {g} \) with its structure
morphisms given by linear maps \( d_{h}:h\rightarrow h[1] \), \( m:\frak {g}\otimes h\rightarrow h \)
one can define \( Q_{\mathfrak {g}}- \)connection \( \mathcal{D}_{Q_{\mathfrak {g}}} \)
on trivial bundle over \( \frak {g}[1] \) with fibers equal to the linear space
\( h \). The covariant derivatives of constant sections are given by 
\begin{equation}
\label{dtwisted}
(\mathcal{D}_{Q_{\mathfrak {g}}}s)(\gamma ):=d_{h}s+m(\gamma \otimes s)\, \, \text {for}\, s\in h
\end{equation}
 Using the Leibnitz rule the covariant derivative \( \mathcal{D}_{Q_{\mathfrak {g}}}s \)
can be defined for arbitrary section of the trivial bundle. The conditions imposed
on linear maps \( (d_{h},m) \) in order that they define dg-module structure
on \( h \) can be reformulated simply as \( [\mathcal{D}_{Q_{\mathfrak {g}}},\mathcal{D}_{Q_{\mathfrak {g}}}]=0 \).
It follows that for \( \gamma \in MC_{\frak {g}}(R) \) one has twisted differential
acting on the fiber over \( \gamma  \) giving rise to complex \( (h\otimes R,\mathcal{D}_{Q_{\mathfrak {g}}}|_{\gamma }) \).

\begin{prop}
\label{H}Given a solution to Maurer Cartan equation \( \gamma \in MC_{\frak {g}}(R) \)
and an element of the group of gauge transformations \( \varphi \in \exp \bigl ((\frak {g}_{A_{\infty }}\otimes \frak {M}_{R})^{0}\bigr ) \)
the cohomology of complexes \( (h\otimes R,\mathcal{D}_{Q_{\mathfrak {g}}}|_{\gamma }) \)
and \( (h\otimes R,\mathcal{D}_{Q_{\mathfrak {g}}}|_{\gamma ^{\varphi }}) \)
are canonically isomorphic and the correspondence which associates to \( \gamma \in MC_{\frak {g}}(R) \)
the cohomology of complex \( (h\otimes R,\mathcal{D}_{Q_{\mathfrak {g}}}|_{\gamma }) \)
defines a sheaf over the moduli space \( \mathcal{M}_{\frak {g}} \). 
\end{prop}
\begin{proof}
Morally this means that \( Q_{\mathfrak {g}}- \)equivariant sheaves descend
naturally to the quotient \( \mathcal{M}_{\frak {g}} \). To check that the
cohomology of complexes \( (h\otimes R,\mathcal{D}_{Q_{\mathfrak {g}}}|_{\gamma }) \)
and \( (h\otimes R,\mathcal{D}_{Q_{\mathfrak {g}}}|_{\gamma ^{\varphi }}) \)
are canonically isomorphic it is sufficient to consider infinitesimal gauge
transformation \( \varphi =\varepsilon \alpha  \), \( \varepsilon ^{2}=0 \).
Then the map \( Id+\varepsilon m(\alpha \bigotimes \cdot ) \) gives the needed
chain isomorphism.
\end{proof}
Same statement holds true for arbitrary \( L_{\infty }- \)modules over \( L_{\infty }- \)algebra
(in this case the vector field \( Q \) and the \( Q- \)connection \( \mathcal{D}_{Q} \)
will be given by formal power series and as usual one should express everything
in terms of co-connections on comodules etc.).

Below we will specify this description in the case of the moduli space \( \mathcal{M} \)
and the sheaf corresponding to the natural extension of the Gauss-Manin local
system.

\subsection{Local system on \protect\( \mathcal{M}^{cs}\protect \). }

Let us give first the description of the Gauss-Manin local system \( (\mathcal{H}^{cs},\nabla ) \)
over the moduli space \( \mathcal{M}^{cs} \). The sheaf \( \mathcal{H}^{cs} \)
is described by \( \frak {g}_{cs}- \)module \( (\Omega ^{*,*},\overline{\partial }+\partial ) \)
with \( \frak {g}_{cs}- \)action \( \gamma \in \frak {g}_{cs} \) \( \rightarrow \mathcal{L}_{\gamma } \),
\( \mathcal{L}_{\gamma }:=[\partial ,i_{\gamma }] \). This description proceeds
as follows. Let \( \gamma \in MC_{cs}(R) \) be a solution to eq.(\ref{mc})
and \( \phi _{\gamma }:(Spec\, R,0)\rightarrow (\mathcal{M}^{cs},[X]) \) be
the corresponding morphism representing a deformation of complex structures
over \( R \). The sheaf \( \mathcal{H}^{cs} \) can be described via its inverse
images \( \phi ^{*}_{\gamma }\mathcal{H}^{cs} \). The sheaf \( \phi ^{*}_{\gamma }\mathcal{H}^{cs}\rightarrow \text {Spec}\, R \)
is given by \( R- \)module defined as the cohomology of the complex 
\begin{equation}
\label{OmR}
(\Omega ^{*,*}_{R},\, d_{\gamma }):=(\Omega ^{*,*}\bigotimes R,\, \, \overline{\partial }+\partial +[i_{\gamma },\partial ])
\end{equation}
 here \( i_{\gamma } \) denotes the operation of contraction with a vector
field extended via wedge multiplication by \( (0,1)- \)form. One can show easily
that 
\[
\phi ^{*}_{\gamma }\mathcal{H}^{cs}=\bigoplus _{w}R^{w}f_{*}(\Omega ^{*}_{\widetilde{X}/\text {Spec}\, R},d_{\widetilde{X}/\text {Spec}\, R})\]
 where \( f:\widetilde{X}\rightarrow \text {Spec}\, R \) is the flat family
of complex structures on \( X \) corresponding to \( \gamma  \). The sheaf
\( \mathcal{H}^{cs} \) is equipped with flat connection 
\begin{equation}
\label{nabla}
\nabla :\mathcal{H}^{cs}\rightarrow \mathcal{H}^{cs}\otimes \Omega ^{1}_{\mathcal{M}^{cs}}\, ,\nabla ^{2}=0
\end{equation}
 Given a solution \( \gamma \in MC_{cs}(R) \) the covariant derivative along
a vector field \( v\in Der(R) \) acting on element \( \varphi  \) of \( \phi ^{*}_{\gamma }\mathcal{H}^{cs} \)
is written as 
\[
\nabla _{v}\varphi =\frac{\partial \varphi }{\partial v}+i_{\frac{\partial \gamma }{\partial v}}(\varphi )\]
 one can check that this is indeed the Gauss-Manin flat connection arising from
the identification of de Rham cohomology with Betti cohomology.

\subsection{Local system on \protect\( \mathcal{M}\protect \).}

Notice that the \( \frak {g}_{cs}- \)action on \( (\Omega ^{*,*},\overline{\partial }+\partial ) \)
extends naturally to the \( \frak {g}_{A_{\infty }}- \)action \( \gamma \rightarrow \mathcal{L}_{\gamma } \),
\( \mathcal{L}_{\gamma }:=[\partial ,i_{\gamma }] \) where \( i_{\gamma } \)
denotes the operator of (holomorphic) contraction \( \Lambda ^{*}T\otimes \Omega ^{*}\rightarrow \Omega ^{*} \)extended
by wedge multiplication by \( (0,q)- \)form. The description of the sheaf \( \mathcal{H} \)
over \( \mathcal{M} \) follows the same scheme as above. Let us notice that
the operator \( d_{\gamma }=\overline{\partial }+\partial +[i_{\gamma },\partial ] \)
has the property \( d_{\gamma }^{2}=0 \) for arbitrary element \( \gamma \in MC_{A_{\infty }}(R) \)
describing an \( A_{\infty }- \)deformation \( \phi _{\gamma }:\text {Spec}\, R\rightarrow \mathcal{M} \).
Therefore one can associate to such element \( \gamma  \) an \( R- \)module
\( ''\phi ^{*}_{\gamma }\mathcal{H}'' \) equal to the cohomology of the complex
\( (\Omega ^{*,*}_{R},\, d_{\gamma }) \). Also given two equivalent solutions
\( \gamma _{1},\gamma _{2} \) the cohomology of the complexes \( (\Omega ^{*,*}_{R},\, d_{\gamma _{1}}),\, (\Omega ^{*,*}_{R},\, d_{\gamma _{2}}) \)
are canonically isomorphic. It follows that the formula (\ref{OmR}) with \( \gamma \in \frak {g}_{A_{\infty }} \)
gives a sheaf \( \mathcal{H}\rightarrow \mathcal{M} \) .

To explain this in more details one has to recall first a definition of a notion
of sheaf over formal moduli space described by functor (\ref{def(R)}). We will
give a sketch of how one should proceed and leave the details to the reader.
First of all let us consider functor \( \widetilde{{Def}}^{\mathbb {Z}}_{\frak {g}} \)
which associates to graded \( \mathbb {C}- \)Artin algebra \( R \) the groupoid(=category
where all morphisms are isomorphisms) whose objects are in one-to-one correspondence
with elements \( \gamma \in MC(R) \) and where morphisms \( \gamma _{1}\rightarrow \gamma _{2} \)
are in one-to-one correspondence with elements \( \widetilde{\alpha }\in \exp \bigl ((\frak {g}\otimes \frak {M}_{R})^{0}\bigr ) \)
transforming \( \gamma _{1} \) to \( \gamma _{2} \). Notice that the functor
\( Def^{\mathbb {Z}}_{\frak {g}} \) can be obtained from \( \widetilde{{Def}}^{\mathbb {Z}}_{\frak {g}} \)
by taking \( \pi _{0} \)=the sets of connected components of objects. By definition
a sheaf over moduli space described by \( Def^{\mathbb {Z}}_{\frak {g}} \)
is a functor which associates to a graded \( \mathbb {C}- \)Artin algebra \( R \)
a representation of the grouppoid \( \widetilde{{Def}}^{\mathbb {Z}}_{\frak {g}}(R) \)
in the category of \( R- \)modules. In our case it is the representation given
by correspondence \( \gamma \rightarrow \phi _{\gamma }^{*}\mathcal{H} \).

The sheaf \( \mathcal{H} \) is equipped with flat connection defined by the
same formula (\ref{nabla}). Again formally this should be understand in terms
of functors on Artin algebras similar to above.

In the sequel we will need a slight generalization of this situation involving
formal parameter \( \hbar  \). The fiber of the sheaf \( \mathcal{H}^{\hbar } \)
over \( \phi _{\gamma }\in \mathcal{M} \) is given by cohomologies of complex
\begin{equation}
\label{drH}
(\Omega ^{*,*}\otimes \mathbb {C}[[\hbar ^{-1},\hbar ]],\overline{\partial }+[\partial ,i_{\gamma }]+\hbar \partial )
\end{equation}
 The covariant derivative along vector field \( v\in Der(R) \) is given by
the formula 
\begin{equation}
\label{hGM}
\nabla _{v}=\frac{\partial }{\partial v}+\frac{1}{\hbar }i_{\frac{\partial \gamma }{\partial v}}
\end{equation}

\subsection{\protect\( \mathcal{D}-\protect \)modules over moduli spaces.}

As a side remark let us mention that more generally, given a moduli space \( \mathcal{M} \)
described by a differential graded Lie algebra \( \frak {g} \), a sheaf \( H \)
equipped with flat connection (\( \mathcal{D}- \)module) over \( \mathcal{M} \)
can be described using the technique of deformation theory as follows. Let us
denote via \( \mathbb {C}[\xi ] \) the differential graded \emph{commutative}
algebra with generator \( \xi ,\, \xi ^{2}=0 \), \( \text {deg}\, \xi =-1 \)
and the differential \( d=\frac{\partial }{\partial \xi } \). Notice that the
space \( \frak {g}\otimes \mathbb {C}[\xi ] \) is naturally a differential
graded Lie algebra.

\begin{prop}
\label{pr:dmod}Given a \( \frak {g}- \)module \( h \) such that \( \frak {g} \)-module
structure on \( h \) is lifted to structure of a module over the differential
graded Lie algebra \( \frak {g}\otimes \mathbb {C}[\xi ] \) one has natural
flat connection \( \nabla  \) on the sheaf \( H \) over \( \mathcal{M}_{\frak {g}} \)associated
with \( h \) according to prop.\ref{H}. The covariant derivative along \( v\in Der(R) \)
acting on cohomology of \( \mathcal{D}_{Q}|_{\gamma }, \) \( \gamma \in MC_{\frak {g}}(R) \)
is given by \( \nabla _{\nu }=\frac{\partial }{\partial v}+i_{\frac{\partial \gamma }{\partial v}} \)where
\( i_{\alpha } \) denotes the action of the element \( \xi \alpha \in \varepsilon \frak {g} \). 
\end{prop}
\begin{proof}
It is straightforward to check that \( [\nabla ,\mathcal{D}_{Q}]=0 \), \( [\nabla ,\nabla ]=0 \).

In our case \( \frak {g}_{A_{\infty }}- \)module is the complex (\ref{drH})
with \( \frak {g}_{A_{\infty }}- \)action \( \gamma \in \frak {g}_{A_{\infty }} \)
\( \rightarrow \mathcal{L}_{\gamma } \), \( \mathcal{L}_{\gamma }:=[\partial ,i_{\gamma }] \).
The operators \( \frac{1}{\hbar }i_{\gamma } \) for \( \varepsilon \gamma \in \varepsilon \frak {g}_{A_{\infty }} \)
define the lifting of this action to \( \frak {g}_{A_{\infty }}\otimes \mathbb {C}[\xi ]- \)action. 
\end{proof}
\begin{rem}
Let us mention that the same result holds true for \( L_{\infty }- \)module
structure over an \( L_{\infty }- \)algebra \( \mathcal{G} \) which can be
lifted to \( L_{\infty }- \)module structure over \( \mathcal{G}\otimes \mathbb {C}[\xi ] \).
The proof is the same with the only modification that \( [\nabla ,\nabla ]=-2[\mathcal{D}_{Q},S] \)
for some operator \( S \). We see that \( L_{\infty }- \)module over \( \mathcal{G}\otimes \mathbb {C}[\xi ] \)
describe higher homotopy generalization of the notion of sheaf with flat connection.
There is in fact an infinite sequence of operators \( \widehat{\nabla }^{0}=\mathcal{D}_{Q}, \)
\( \widehat{\nabla }^{1}=\nabla  \), \( \widehat{\nabla }^{2}=S \), \( \widehat{\nabla }^{3} \)\( \ldots  \)
satisfying higher homotopy identities. The existence of such higher homotopies
allows one to define canonically a parallel transport operator acting on the
complex \( (h,D_{Q}) \) for the connection which is flat only up to homotopy.
The details will appear in \cite{B3}. 
\end{rem}

\section{Semi-infinite Hodge structure\label{sect:semiHodge}}

In this section we define an analog of Hodge filtration for \( A_{\infty }- \)deformations
\( \phi \in \mathcal{M} \).

\subsection{Grassmanian\label{grassm}}

For a \( \mathbb {Z}_{2}- \)graded \( \mathbb {C}- \)vector space \( H=(H^{\text {even}},H^{\text {odd}}) \),
denote via \( H[[\hbar ^{-1},\hbar ]] \) the space of formal series 
\[
\sum ^{+\infty }_{i=-\infty }v_{i}\hbar ^{i}+\sum _{j=-\infty }^{+\infty }u_{j}\hbar ^{j+\frac{1}{2}},\, v_{i}\in H^{\text {even}},u_{j}\in H^{\text {odd}}\]
 Denote also 
\[
H[[\hbar ]]=\{\sum _{r\geq 0}a_{r}\hbar ^{r}\},\: H[[\hbar ^{-1}]]=\{\sum _{r<0}a_{r}\hbar ^{r}\}\subset H[[\hbar ^{-1},\hbar ]]\]
 and let \( H((\hbar )) \) be the analogous ring of Laurent power series and
\( p_{+} \) be the linear projection \( H[[\hbar ^{-1},\hbar ]]\rightarrow H[[\hbar ]] \)
along \( H[[\hbar ^{-1}]] \).

Let \( \text {dim}_{\mathbb {C}}H<\infty  \). Let us introduce grassmanian
\( Gr_{\frac{\infty }{2}}(H[[\hbar ^{-1},\hbar ]]) \) of semi-infinite subspaces
in \( H[[\hbar ^{-1},\hbar ]] \). Firstly it contains grassmanian \( Gr_{\frac{\infty }{2}}\left( H((\hbar ))\right)  \)
which consists of semi-infinite subspaces in \( H((\hbar )) \). In fact, as
we work in the neighborhood of moduli space of complex structures, we will be
dealing only with the points of \( Gr_{\frac{\infty }{2}}(H[[\hbar ^{-1},\hbar ]]) \)
which lie in the formal neighborhood\footnote{%
In general, for deformations which are not of finite order the subspaces which
we study do not belong to \( Gr_{\frac{\infty }{2}}\left( H((\hbar ))\right)  \).
One can repeat the same story in the analytic setting i.e. on the level of convergent
solutions to Maurer-Cartan equation. In this case some type of Segal-Wilson
grassmanian should be used. Then subspaces associated analogously with such
convergent non-commutative deformations (see below) will not in general lie
in \( H((\hbar )) \) as well. 
} of \( Gr_{\frac{\infty }{2}}\left( H((\hbar ))\right)  \). Under points of
\( Gr_{\frac{\infty }{2}}(H[[\hbar ^{-1},\hbar ]]) \) over (pro-)Artin algebra
\( R \) we will always mean in this paper points of such neighborhood. More
precisely, a point \( L\in Gr_{\frac{\infty }{2}}(H[[\hbar ^{-1},\hbar ]])(R) \)
over an Artin algebra \( R \) (i.e. a \char`\"{}family of subspaces\char`\"{}
represented by morphism \( Spec\, R\rightarrow Gr_{\frac{\infty }{2}}(H[[\hbar ^{-1},\hbar ]]) \))
is an \( R\otimes \mathbb {C}[[\hbar ]]- \)submodule \( L\subset H((\hbar ))\otimes R \)
such that 
\[
L\otimes _{\mathbb {C}[[\hbar ]]}\mathbb {C}((\hbar ))=H((\hbar ))\otimes R\]
 A point \( L\in Gr_{\frac{\infty }{2}}(H[[\hbar ^{-1},\hbar ]])(R) \) over
a pro-Artin algebra \( R \) is a submodule in \( H[[\hbar ^{-1},\hbar ]]\widehat{\otimes }R \)
given by projective limit of submodules \( (proj)\lim _{i}\, L\otimes _{R}R/\mathfrak {M}_{R}^{i}\in Gr_{\frac{\infty }{2}}(H[[\hbar ^{-1},\hbar ]])(R/\mathfrak {M}_{R}^{i}) \)
defined before. Such submodules always satisfy conditions \( \dim _{R}\ker (p_{+}|_{L})<\infty  \),
\( \dim _{R}\, \text {coker}\, (p_{+}|_{L})<\infty  \).

Let \( F \) be a decreasing filtration on \( H \), \( \text {dim}_{\mathbb {C}}H<\infty  \).
A filtration in the category of \( \mathbb {Z}_{2}- \)graded vector spaces
is given by a pair of filtrations \( F=(F^{even},F^{odd}) \) on \( H^{even} \)and
\( H^{odd} \) respectively. Assume in addition that \( F^{even} \) has indexes
taking values in \( \mathbb {Z} \) and \( F^{odd} \) has indexes taking values
in \( \frac{1}{2}+\mathbb {Z} \). To shorten the notations \( F^{\geq r} \)will
denote below the subspace \( (F^{even})^{r} \) if \( r\in \mathbb {Z} \) and
the subspace \( (F^{odd})^{r} \) if \( r\in \frac{1}{2}+\mathbb {Z} \). To
such filtration one can associate a subspace \( L^{F}\in Gr_{\frac{\infty }{2}}(H[[\hbar ^{-1},\hbar ]]) \):
\[
L^{F}:=\text {linear\, span}_{r\in \mathbb {Z}[\frac{1}{2}]}F^{\geq r}\hbar ^{-r}[[\hbar ]]\]
 The correspondence \( F\rightarrow L^{F} \) is injective. In particular one
has the subspace which corresponds to Hodge filtration\footnote{%
we use indexes which are shifted with respect to the standard notations 
} on \( H^{*}(X,\mathbb {C}) \)

\begin{equation}
\label{Hfiltr}
F^{\geq r}(\phi ):=\bigoplus _{p-q\geq 2r,\, p-q\equiv 2r}H^{p,q},\: H^{p,q}=\{[\varphi ]\in H^{p+q}(X,\mathbb {C})[-p-q]|\varphi \in \Omega ^{p,q}(\phi )\}
\end{equation}
 which is associated with a deformation of complex structure \( \phi \in \mathcal{M}^{cs} \).

\subsection{Semi-infinite Hodge structure. }

We claim that there is a subspace \( L(\phi )\in Gr_{\frac{\infty }{2}}(H^{*}(X,\mathbb {C})[[\hbar ^{-1},\hbar ]])(R) \)
associated canonically to any \( A_{\infty }- \)deformation \( \phi \in \mathcal{M}(R) \).
For \( \phi \in \mathcal{M}^{cs}(R) \) it coincides with the subspace \( L^{F(\phi )} \)corresponding
to the Hodge filtration (\ref{Hfiltr}).

The subspace \( L(\phi ) \) is defined as follows. Let \( \gamma \in MC(R) \)
describes an \( A_{\infty }- \)defor\-mation \( \phi \in \mathcal{M}(R) \)
where \( R \) is an Artin algebra or a projective limit of Artin algebras.
To such deformation we associate an \( R- \)module \( L(\phi ) \) in \( H^{*}(X,\mathbb {C})[[\hbar ^{-1},\hbar ]]\otimes R \)
representing morphism \( \text {Spec}\, R\rightarrow Gr_{\frac{\infty }{2}}(H^{*}(X,\mathbb {C})[[\hbar ^{-1},\hbar ]]) \).
Let \( l_{\hbar } \) be rescaling operator: 
\[
l_{\hbar }:\varphi ^{p,q}\rightarrow \hbar ^{\frac{q-p}{2}}\varphi ^{p,q}\]
 Denote \( \widetilde{L}_{\gamma }\subset \Omega ^{*}(X,\mathbb {C})[[\hbar ^{-1},\hbar ]]\bigotimes R \)
the subspace of elements of the form\footnote{%
The action of \( \exp (\frac{1}{\hbar }i_{\gamma }) \) on arbitrary element
from \( \Omega ^{*}(X,\mathbb {C})\widehat{\bigotimes }k[[\hbar ]] \) is well-defined
even in the case when \( R \) is a \textsl{pro}-Artin algebra since \( \exp (\frac{1}{\hbar }\gamma ) \)
belongs to subspace of elements of the form \( \sum _{j\geq 0}\gamma _{j}\, \, \gamma _{j}\in \mathfrak {M}^{j}_{R}\bigotimes \hbar ^{-j}\frak {g}_{A_{\infty }}\widehat{\otimes }k[[\hbar ]] \). 
}
\begin{equation}
\label{Lgamma}
l_{\hbar }\exp (\frac{1}{\hbar }i_{\gamma })(\varphi _{0}+\varphi _{1}\hbar +\varphi _{2}\hbar ^{2}\ldots )\: \: \text {where}\, \varphi _{i}\in \Omega ^{*}(X,\mathbb {C})\otimes R
\end{equation}

\begin{prop}
The subspace \( \widetilde{L}_{\gamma } \) is a sub-complex of the complex
\begin{equation}
\label{dRh}
(\Omega ^{*}(X,\mathbb {C})[[\hbar ^{-1},\hbar ]]\otimes R,\, \hbar ^{\frac{1}{2}}(\overline{\partial }+\partial ))
\end{equation}

\end{prop}
\begin{proof}
\begin{eqnarray}
l_{\hbar }^{-1}\hbar ^{\frac{1}{2}}(\overline{\partial }+\partial )l_{\hbar } & = & (\overline{\partial }+\hbar \partial )\\
\exp (-\frac{1}{\hbar }i_{\gamma })(\overline{\partial }+\hbar \partial )\exp (\frac{1}{\hbar }i_{\gamma }) & = & \overline{\partial }+\hbar \partial +\frac{1}{\hbar }(i_{\overline{\partial }\gamma +\frac{1}{2}[\gamma ,\gamma ]})+[\partial ,i_{\gamma }]=\label{e(-gamma)de(gamma)} \\
 & = & \overline{\partial }+[\partial ,i_{\gamma }]+\hbar \partial \nonumber 
\end{eqnarray}

\end{proof}
Notice that the cohomology groups of the complex (\ref{dRh}) are isomorphic
to \break \( H^{*}(X,\mathbb {C})[[\hbar ^{-1},\hbar ]]\bigotimes R \).

\begin{defn*}
The \( R- \)submodule \( L(\phi )\subset H^{*}(X,\mathbb {C})[[\hbar ^{-1},\hbar ]]\bigotimes R \)
consists of cohomology classes represented by elements from \( \widetilde{L}_{\gamma } \). 
\end{defn*}
\begin{thm}
To an arbitrary \( A_{\infty }- \)deformation \( \phi \in \mathcal{M}(R) \)
there is canonically associated subspace \( L(\phi )\in Gr_{\frac{\infty }{2}}(H^{*}(X,\mathbb {C})[[\hbar ^{-1},\hbar ]])(R) \)
defined by formula (\ref{Lgamma}). For \( \phi \in \mathcal{M}^{cs}(R) \)
this is the subspace \( L^{F(\phi )} \)corresponding to the Hodge filtration
(\ref{Hfiltr}). 
\end{thm}
\begin{proof}
We need to check first that \( L(\phi )=L(\phi ') \) for \( \gamma \sim \gamma ' \).
It is enough to prove this for \( \gamma '=\gamma +\varepsilon (\overline{\partial }\alpha +[\gamma ,\alpha ]),\, \varepsilon ^{2}=0 \).
Applying the standard commutation rules 
\[
[\mathcal{L}_{\gamma _{1}},i_{\gamma _{2}}]=i_{[\gamma _{1},\gamma _{2}]},\, \, \, [i_{\gamma _{1}},i_{\gamma _{2}}]=0\]
 one gets: 
\begin{equation}
\label{dexp}
\exp (\frac{1}{\hbar }i_{\gamma +\varepsilon \overline{\partial }\alpha +[\gamma ,\alpha ]})(\text {Id}+\varepsilon [\partial ,i_{\alpha }])=\exp (\frac{1}{\hbar }i_{\gamma })(\text {Id}+\varepsilon [\overline{\partial }+[\partial ,i_{\gamma }]+\hbar \partial ,\frac{1}{\hbar }i_{\alpha }])
\end{equation}
 Assume that \( \varphi \in \text {ker}(\overline{\partial }+\partial )\cap \widetilde{L}_{\gamma } \),
or equivalently 
\[
\exp (-\frac{1}{\hbar }i_{\gamma })l_{\hbar }^{-1}\varphi \in \text {ker}(\overline{\partial }+[\partial ,i_{\gamma }]+\hbar \partial )\cap \Omega ^{*}(X,\mathbb {C})\widehat{\otimes }\mathbb {C}[[\hbar ]]\]
 Applying the operator (\ref{dexp}) to \( \widetilde{\varphi }=\exp (-\frac{1}{\hbar }i_{\gamma })l_{\hbar }^{-1}\varphi  \)
we get 
\[
\varphi \pm \varepsilon l_{\hbar }\exp (\frac{1}{\hbar }i_{\gamma })(\overline{\partial }+[\partial ,i_{\gamma }]+\hbar \partial )\frac{1}{\hbar }i_{\alpha }\widetilde{\varphi }\in \widetilde{L}_{\gamma +\varepsilon (\overline{\partial }\alpha +[\gamma ,\alpha ])}\]
 on the other hand 
\[
l_{\hbar }\exp (\frac{1}{\hbar }i_{\gamma })(\overline{\partial }+[\partial ,i_{\gamma }]+\hbar \partial )=(\overline{\partial }+\partial )\hbar ^{\frac{1}{2}}l_{\hbar }\exp (\frac{1}{\hbar }i_{\gamma })\]
 Let \( \gamma \in \Omega ^{0,1}(X,T)\otimes \mathfrak {M}_{R^{0}} \) be a
solution to Maurer-Cartan equation describing a deformation of complex structure
\( \phi :\text {Spec}\, R\rightarrow \mathcal{M}^{cs} \). Let \( \{z^{i}\} \)
be a set of complex coordinates with respect to the fixed initial complex structure
\( J_{0} \). The differentials of the deformed complex coordinates can be written
as \( (dz^{i})^{\text {new}}=dz^{i}+\sum _{\overline{j}}\gamma _{\overline{j}}^{i}d\overline{z}^{\overline{j}} \),
where \( \gamma =\sum _{i,\overline{j}}\gamma ^{i}_{\overline{j}}d\overline{z}^{\overline{j}}\frac{\partial }{\partial z^{i}} \).
Therefore \( F^{\geq r}(\phi )=\exp i_{\gamma _{*}}F^{\geq r}(J_{0}) \). Notice
also that for \( \gamma \in \Omega ^{0,1}(X,T)\otimes \mathfrak {M}_{R^{0}} \)
one has \( l_{\hbar }\circ \exp (\frac{1}{\hbar }i_{\gamma })=\exp (i_{\gamma })\circ l_{\hbar } \).
It follows that \( L(\phi )=L^{F(\phi )} \) for \( \phi \in \mathcal{M}^{cs}(R) \).

It is easy to see that there exist elements \( \alpha ^{k}_{j}\in \Omega ^{*}(X) \),
\( k=1,\ldots ,\text {dim}\, H^{*}(X,\mathbb {C}) \), \( j=0,\ldots ,m(k) \),
such that \( \overline{\partial }\alpha ^{k}_{j}=\partial \alpha ^{k}_{j-1} \),
\( \overline{\partial }\alpha ^{k}_{0}=0 \), \( \partial \alpha ^{k}_{m(k)}=0 \)
and such that \break \( \{\sum ^{m(k)}_{j=0}\alpha ^{k}_{j}\}_{k=1,\ldots ,\text {dim}\, H^{*}(X,\mathbb {C})} \)
give a basis in \( H^{*}(X,\mathbb {C}) \). Then given a solution to Maurer-Cartan
equation \( \gamma \in MC_{A_{\infty }}(R) \) over an Artin algebra \( R \)
describing \( \phi _{\gamma }\in \mathcal{M}(R) \) one has 
\[
\hbar ^{N}\exp (-\frac{1}{\hbar }i_{\gamma })\sum _{j=0}^{m(k)}\alpha ^{k}_{j}\hbar ^{j}\in \Omega ^{*}(X){\otimes }\mathbb {C}[[\hbar ]]\]
 for sufficiently large \( N \). Therefore \( l_{\hbar }\sum _{j=0}^{m(k)}\alpha ^{k}_{j}\hbar ^{j}\in \hbar ^{-N}\widetilde{L}_{\gamma } \)
and \( L(\phi _{\gamma })\otimes _{\mathbb {C}[[\hbar ]]}\mathbb {C}((\hbar ))=H^{*}(X,\mathbb {C})((\hbar ))\otimes R \)
\end{proof}
\begin{rem}
The operator \( \frac{1}{\hbar }i_{\gamma } \) is an operator of the flat connection
on the sheaf of cohomologies of complexes (\ref{drH}) described in the previous
section. The operator \( l_{\hbar } \) identifies cohomology of complexes (\ref{drH})
and (\ref{dRh}). We see that \( L(\phi ) \) is the subspace obtained as a
result of parallel transport to the base point \( [X]\in \mathcal{M} \) of
the subspace \( \{[\varphi _{0}+\varphi _{1}\hbar +\varphi _{2}\hbar ^{2}\ldots ]\} \)
living over \( \phi \in \mathcal{M} \) and identification provided by \( l_{\hbar } \). 
\end{rem}
\begin{rem} If \( X \) is a compact Kahler manifold then \( \overline{\partial }+[\partial ,i_{\gamma }] \)
and \( \partial  \) satisfy \( \partial \overline{\partial }- \)lemma. It
follows that \( L(\phi ) \) is a free \( R- \)module for any \( \phi _{[\gamma ]}\in \mathcal{M}(R) \).\end{rem}

\begin{rem}As we already mentioned the moduli space \( \mathcal{M} \) can be
viewed as a moduli space parameterizing \( A_{\infty }- \)deformations of \( D^{b}Coh(X) \).
The local system over \( \mathcal{M} \) is the local system of periodic cyclic
homology of the \( A_{\infty }- \)categories and the semi-infinite subspace
coincides with negative cyclic homology. Operator \( l_{\hbar } \) comes from
natural grading on periodic cyclic homology of \( D^{b}Coh(X) \). \end{rem}

\section{Griffiths transversality \label{section griffiths}}

In this section we study variations of semi-infinite subspaces introduced in
the previous section. Recall that the family of filtrations \( F^{\geq r}(\phi ) \)
for deformations of complex structure \( \phi \in \mathcal{M}^{cs} \) satisfies
Griffiths transversality condition with respect to the Gauss-Manin connection:
\[
\nabla F^{\geq r}\subset F^{\geq r-1}\]

\begin{thm}
\label{th:griff}The covariant derivative of \( L(\phi )\in Gr_{\frac{\infty }{2}}(H^{*}(X,\mathbb {C})[[\hbar ^{-1},\hbar ]])(R),\, \phi \in \mathcal{M}(R) \)
with respect to the Gauss-Manin connection satisfies 
\begin{equation}
\label{eq:grtransv}
\nabla L\subset \hbar ^{-1}L
\end{equation}

\end{thm}
\begin{proof}
Let \( \gamma \in MC(R) \) describes an \( A_{\infty }- \)deformation \( \phi _{\gamma }\in \mathcal{M}(R) \).
Let 
\[
[\varphi ]=[l_{\hbar }\exp (\frac{1}{\hbar }i_{\gamma })(\varphi _{0}+\varphi _{1}\hbar +\ldots )]\in L(\phi _{\gamma })\subset H^{*}(X,\mathbb {C})[[\hbar ^{-1},\hbar ]]\otimes R\]
 describes a family of elements of the varying semi-infinite subspaces. Then
for any vector field \( v\in Der(R) \): \begin{multline}
\nabla _{v}[\varphi ]=  [l_{\hbar }\exp (\frac{1}{\hbar }i_{\gamma
}) (\partial _{v}\varphi _{0}+\partial _{v}\varphi _{1}\hbar
+\ldots ) +\\
+\frac{1}{\hbar }l_{\hbar }\exp (\frac{1}{\hbar
}i_{\gamma }) (i_{\partial _{v}\gamma }\varphi _{0}+ i_{\partial
_{v}\gamma }\varphi _{1}\hbar +\ldots )]\in \hbar ^{-1}L(\phi_{\gamma})
\end{multline}
\end{proof}
Notice that one has induced map (symbol of the Gauss-Manin connection) 
\begin{equation}
\label{symbolGMgen}
T_{\phi }\mathcal{M}\otimes L(\phi )/\hbar L(\phi )\rightarrow \hbar ^{-1}L(\phi )/L(\phi )
\end{equation}
 In the case of deformations of complex structure \( \phi \in \mathcal{M}^{cs} \)
this is the standard map 
\begin{equation}
\label{symbolGM}
\left( \bigoplus _{p,q}H^{q}(X,\Lambda ^{p}T_{X})[2]\right) \otimes \left( \bigoplus _{r}F^{\geq r}/F^{\geq r-1}\right) \rightarrow \left( \bigoplus _{r}F^{\geq r}/F^{\geq r-1}\right) 
\end{equation}

\begin{rem}Let us consider integrals 
\begin{equation}
\label{periods}
\varphi ^{i}(t,\hbar )=\int _{\Delta _{i}}\varphi (t,\hbar )\, \, \, \text {for}\, \, \, \varphi (t,\hbar )\in L(t),
\end{equation}
 where \( \{\Delta _{i}\} \) is a basis in \( H^{*}(X,\mathbb {Z}) \) and
\( t\in \mathcal{M}(R) \) represents a class of \( A_{\infty }- \)deforma\-tions
of complex projective manifold \( X \). Such integrals satisfy system of differential
equations which are generalizations to the case of non-commutative deformations
of Picard-Fuchs equations for usual periods. This can be seen as follows. Assume
that we are given elements \( \varphi _{\alpha }(0)\in F^{\geq r_{\alpha }}(0), \)\( \varphi _{\alpha }(0)\notin F^{r_{\alpha }-1}(0) \)
which form a basis in \( H^{*}(X,\mathbb {C}) \). Let \( W_{\leq r} \) be
an increasing \( \mathbb {Z}[\frac{1}{2}]- \)graded filtration opposite to
filtration \( F^{\geq r}(0) \): that is, for any \( r \), \( \bigoplus _{i\equiv 2r}H^{i}(X,\mathbb {C})=F^{\geq r}\oplus W_{\leq r} \).
Then \( L_{W}:=\text {linear\, span}_{r\in \mathbb {Z}[\frac{1}{2}]}W_{\leq r}\hbar ^{-r}[[\hbar ^{-1}]] \)
is a subspace transversal to \( L^{F(0)} \). Let us consider elements \( \widetilde{\varphi }_{\alpha }(t,\hbar )=L(t)\cap (\varphi _{\alpha }(0)\hbar ^{-r_{\alpha }}+L_{W}) \).
Notice that for \( t\in \mathcal{M}^{cs} \), \( \widetilde{\varphi }_{\alpha }(t,\hbar )=\hbar ^{-r_{\alpha }}\varphi _{\alpha }(t) \),
where \( \varphi _{\alpha }(t)=F^{\geq r_{\alpha }}(t)\cap (\varphi _{\alpha }(0)+W_{r_{\alpha }}) \).
Then one can show (see proof of prop.\ref{Pr_pic-fuchs} below for analogous
arguments) that the periods \( \widetilde{\varphi }^{i}_{\alpha }(t,\hbar ) \)
satisfy the following system of equations: for any vector field \( v\in Der(R) \)
one has 
\[
\frac{\partial \widetilde{\varphi }_{\alpha }^{i}}{\partial v}=\hbar ^{-1}\sum _{\beta }\Gamma ^{\beta }_{v}(t)\widetilde{\varphi }_{\beta }^{i}\]
 where \( \hbar ^{-1}\Gamma _{v}(t) \) is the \( 1- \)form representing Gauss-Manin
connection (\ref{hGM}).\end{rem}

\begin{rem}
One can generalize easily the notion of real polarized variation of Hodge structure
to the semi-infinite context. For example, the semi-infinite analog of the standard
condition\footnote{%
here for reader convinience we use the standard grading on components of Hodge
filtration 
} 
\[
\forall \, p,k\, \, \, \, \, F^{p}\oplus \overline{F}^{k-p}=H^{k}\]

\end{rem}
is given by 
\[
L(\hbar )\oplus \overline{L}(\overline{\hbar }|_{\overline{\hbar }=\hbar ^{-1}})=H^{*}(X,\mathbb {C})[[\hbar ^{-1},\hbar ]]\]
 satisfied by subspaces \( L(t) \), \( t\in \mathcal{M} \). Similarly, one
can define semi-infinite analog of polarization form.

\section{Quantum periods and counting of rational Gromov-Witten invariants on Calabi-Yau
manifolds. \label{mirror}}

In three dimensions mirror symmetry predictions express rational Gromov-Witten
invariants of a Calabi-Yau threefold in terms of variations of Hodge structure
associated with the deformations of complex structure on the mirror dual Calabi-Yau
threefold (see \cite{CGOP}). We explain in this section that the ``non-commutative''
variations of semi-infinite Hodge structures described above play the same role
in higher dimensional mirror symmetry. In this section \( X \) denotes compact
Kahler manifold with \( c_{1}(T_{X})\in Pic(X) \) equals to zero.

\subsection{Moduli spaces of \protect\( A_{\infty }-\protect \)deformations of Calabi-Yau
manifolds}

In the case when \( c_{1}(T_{X})\in Pic(X) \) equals to zero the moduli space
\( \mathcal{M} \) has especially nice properties. Let us denote 
\[
\mathbf{H}=\bigoplus _{p,q}\mathbb {C}^{dim_{\mathbb {C}}H^{q}(X,\Lambda ^{p}T)}[p-q]\]
 Denote also via \( \mathbb {C}[[t_{\bold {H}}]] \) the algebra of formal power
series on the graded vector space \( \bold {H} \). It is convenient to fix
some choice of a set \( \{t^{a}\} \) of linear coordinates on \( \bold {H} \).

\begin{prop}
\label{M_CY}(\cite{BK}, lemma 2.1) The deformation functor \( {Def}^{\mathbb {Z}}_{\frak {g}_{A_{\infty }}} \)describing
the \( A_{\infty }- \)deformations of \( X \) is isomorphic to the functor
represented by the pro-Artin algebra \( \mathbb {C}[[t_{\bold {H}}]] \). Equivalently,
there exists mini-versal solution to the Maurer-Cartan equation 
\[
\widehat{\gamma }(t)=\sum _{a}\widehat{\gamma }_{a}t^{a}+\frac{1}{2!}\sum _{a_{1},a_{2}}\widehat{\gamma }_{a_{1}a_{2}}t^{a_{1}}t^{a_{2}}+\ldots \, \, \in MC(\mathbb {C}[[t_{\bold {H}}]])\]
 and the cohomology classes \( [\widehat{\gamma }_{a}] \) form a basis in the
cohomology of the complex \( (\frak {g}_{A_{\infty }},\overline{\partial }) \). 
\end{prop}
\begin{rem}
The reader shoud be warned that such solution is by no means unique. A choice
of class of gauge equivalence of such solutions is equivalent to choice of system
of local coordinates: \( (\bold {H},0)\rightarrow (\mathcal{M},[X]) \)
\end{rem}

\subsection{Period map on \protect\( \mathcal{M}^{cs}\protect \).}

The condition \( c_{1}(T_{X})=0 \) implies that there exists ``holomorphic
volume form'' \( \Omega \in \Gamma (X,\Omega ^{n}) \) , \( n=\text {dim}_{\mathbb {C}}X \),
which is non-vanishing at any point of \( X \). Such form is defined up to
multiplication by a non-zero constant. Assume now that we have a family \( \widetilde{X}_{t} \)
of deformations of complex structure on \( \widetilde{X}_{0}=X \). If we choose
a holomorphic volume form \( \Omega _{0} \) on \( \widetilde{X}_{0} \) and
a hyperplane \( l\subset H^{n}(X,\mathbb {C}) \) which is transversal to the
line \( \mathbb {C}\cdot \Omega _{0} \) then the above mentioned ambiguity
in the choice of holomorphic volume element can be fixed using the condition
\[
\Omega ^{l}_{t}-\Omega _{0}\in l\]
 for all values of the parameter \( t \) for which the hyperplane \( l \)
rests transversal to the line \( \mathbb {C}\cdot \Omega _{t} \). The correspondence
\( [\widetilde{X}_{t}]\rightarrow [\Omega ^{l}_{t}] \) defines the period map:
\( \mathcal{M}^{cs}\rightarrow H^{n}(X,\mathbb {C}) \).

\subsection{Period map on \protect\( \mathcal{M}\protect \)\label{permapM}.}

In order to define analogous map on the moduli space \( \mathcal{M} \) let
us fix a choice of subspace in \( H^{*}(X,\mathbb {C})[[\hbar ^{-1},\hbar ]] \)
transversal to \( L^{F(0)} \) where \( F(0) \) is the initial Hodge filtration
on \( H^{*}(X,\mathbb {C}) \). Such subspace can be associated naturally with
an increasing \( \mathbb {Z}[\frac{1}{2}]- \)graded filtration \( W=(W_{even},W_{odd}) \)
on \( H^{*}(X,\mathbb {C}) \) which is opposite to \( F(0) \) in the following
sense: 
\begin{equation}
\label{oppos}
\forall \, r:\, \, \, \bigoplus _{i\equiv 2r}H^{i}(X,\mathbb {C})=F^{\geq r}\oplus W_{\leq r}
\end{equation}
 In the applications to mirror symmetry \( W \) will be the limiting weight
filtration associated with a maximal unipotency cusp on moduli space of complex
structures on \( X \). The subspace \( L_{W}\subset H^{*}(X,\mathbb {C})[[\hbar ^{-1},\hbar ]] \)
associated with the filtration \( W \) is defined in a similar manner as above
for the decreasing filtration \( F^{\geq r} \) (see section \ref{grassm})
\[
L_{W}:=\text {linear\, span}_{r\in \mathbb {Z}[\frac{1}{2}]}W_{\leq r}\hbar ^{-r}[[\hbar ^{-1}]]\]
 where we assumed as above that \( W_{\leq r}:=W^{even}_{\leq r} \) for \( r\in \mathbb {Z} \),
and \( W_{\leq r}:=W^{odd}_{\leq r} \) for \( r\in \frac{1}{2}+\mathbb {Z} \).
It follows from the condition (\ref{oppos}) that \( L_{W} \) is transversal
to \( L^{F(0)} \): 
\begin{equation}
\label{transv}
H^{*}(X,\mathbb {C})[[\hbar ^{-1},\hbar ]]=L_{W}\oplus L^{F(0)}
\end{equation}

Assume now that in addition to filtration \( W \) satisfying (\ref{oppos})
a holomorphic volume element \( \Omega _{0} \) for initial complex structure
on \( X \) is fixed. Let \( \gamma \in MC_{A_{\infty }}(R) \) represents an
\( A_{\infty }- \)deformation \( \phi \in :\mathcal{M}(R). \) Recall that
we denoted via \( L(\phi )\in H^{*}(X,\mathbb {C})[[\hbar ^{-1},\hbar ]]\otimes R \)
an \( R- \)submodule representing the family of semi-infinite subspaces.

\begin{prop}
The intersection \( L(\phi )\bigcap (\Omega _{0}\hbar ^{-\frac{n}{2}}+L_{W}\otimes R) \)
consists of a single element. 
\end{prop}
\begin{proof}
It follows from the condition (\ref{oppos}) that \( L_{W} \) is transversal
to \( L^{F(0)} \). One can show that \( R- \)module \( L(\phi ) \) is free
(it follows, for example, from \cite{B2}, prop. 4.2.1). Therefore \( H^{*}(X,\mathbb {C})[[\hbar ^{-1},\hbar ]]\otimes R=L(\phi )\bigoplus (L_{W}\otimes R) \). 
\end{proof}
We will denote the element defined by intersection \( L(\phi _{\gamma })\bigcap (\Omega _{0}\hbar ^{-\frac{n}{2}}+L_{W}\otimes R) \)
via \( \Psi ^{W}(\phi ,\hbar ) \). For a mini-versal solution \( \widehat{\gamma }\in MC(\mathbb {C}[t_{\bold {H}}]) \)
this element coincides with the series \( (\Phi ^{W}(t_{\bold {H}},\hbar )\vdash \Omega _{0}+\Omega _{0})\hbar ^{-\frac{n}{2}} \)
from \cite{B2} (here \( \Phi ^{W}(t_{\bold {H}},\hbar ) \) is the series defined
in eq. (4.28) after theorem 1 in \cite{B2}). The corresponding morphism 
\[
\Psi ^{W}(\phi _{\widehat{\gamma }},\hbar ):\mathcal{M}\rightarrow H^{*}(X,\mathbb {C})[[\hbar ^{-1},\hbar ]][n(\text {mod})2]\]
 does not depend on the choice of the mini-versal solution \( \widehat{\gamma } \).

\begin{prop}
\label{Pr:PiW}(\cite{B2} propositions 4.2.4 and 4.2.6) The map \( \Psi ^{W}(\phi _{\widehat{\gamma }},\hbar ) \)
has the following properties: 
\end{prop}
\begin{lyxlist}{00.00.0000}
\item [a)]the restriction \( \Pi ^{W}:= \) \( \Psi ^{W}(\hbar )|_{\hbar =1} \)is
well-defined and gives a local isomorphism of formal super-manifolds (local
Torelli theorem): 
\[
\Pi ^{W}:(\mathcal{M},[X])\simeq (H^{*}(X,\mathbb {C})[n(\text {mod})2],[\Omega _{0}])\]

\item [b)]the restriction \( \Pi ^{W}|_{\mathcal{M}^{cs}} \) coincides with the classical
period map \( [X_{t}]\rightarrow [\Omega _{[X_{t}]}^{W_{\leq \frac{n}{2}}}]\in H^{n}(X,\mathbb {C}) \)
defined with help of the hyperplane \( W_{\leq \frac{n}{2}}\bigcap H^{n}(X,\mathbb {C}) \). 
\end{lyxlist}
Let us consider ``value at \( \hbar =\infty  \)'' of \( \Psi ^{W} \) which
is a morphism \( \Psi ^{W}(\infty ):\mathcal{M}\rightarrow Gr\, W \) into associated
quotient of the filtration \( W \) defined as the composition of \( \Psi ^{W} \)
with affine map \( \Omega _{0}\hbar ^{-\frac{n}{2}}+L_{W}\rightarrow L_{W}/\hbar ^{-1}L_{W}=\bigoplus _{r}W_{\leq r}/W_{\leq r-1} \)
sending \( v=\Omega _{0}\hbar ^{-\frac{n}{2}}+v_{-r}\hbar ^{-r}+v_{-r+1}\hbar ^{-r+1}+\ldots \in \Omega _{0}\hbar ^{-\frac{n}{2}}+W_{\leq r}\hbar ^{-r}[[\hbar ^{-1}]] \)
to \( [v_{-r}]\in W_{\leq r}/W_{\leq r-1} \). This morphism is a local isomorphism
as well (\cite{B2}, eq.(4.25)) and defines a set of local coordinates on \( \mathcal{M} \)
which we denote by \( t_{W} \).

\subsection{Mirror symmetry in higher dimensions.}

The map 
\[
\left( \bigoplus _{p,q}H^{q}(X,\Lambda ^{p}T_{X})[-p-q(\text {mod})2]\right) \otimes [\Omega _{0}]\rightarrow \left( \bigoplus _{r}F^{\geq r}/F^{\geq r-1}\right) [n(\text {mod})2]\]
 given by restriction of the symbol of Gauss-Manin connection (\ref{symbolGM})
is an isomorphism. Therefore the set \( \{[\frac{\partial \Psi ^{W}(t_{\bold {H}},\hbar )}{\partial t_{\bold {H}}^{a}}]\, \text {mod}\, L([\widehat{\gamma }])\in \hbar ^{-1}L([\widehat{\gamma }])/L([\widehat{\gamma }])\} \)
is a set of free generators of \( \left( \hbar ^{-1}L([\widehat{\gamma }])\right) /L([\widehat{\gamma }]) \).
This implies that the quantum periods given by components of the map \( \Psi ^{W} \)
satisfy a system of differential equations (non-commutative extension of a version
of Picard-Fuchs equations):

\begin{prop}
\label{Pr_pic-fuchs}The series \( \Psi ^{W}(t_{W},\hbar ) \) satisfy: 
\end{prop}
\[
\frac{\partial ^{2}\Psi ^{W}}{\partial t^{a}_{W}\partial t^{b}_{W}}=\hbar ^{-1}\sum _{c}A_{ab}^{c}(t_{W})\frac{\partial \Psi ^{W}}{\partial t_{W}^{c}}\]
 for some \( A_{ab}^{c}(t_{W})\in \mathbb {C}[[t_{W}]] \) where \( \{t_{W}\} \)
is a set of linear coordinates on \( Gr\, W \) induced by the map \( \Psi ^{W}(\infty ) \).

\begin{proof}
It follows from theorem \ref{th:griff} that 
\[
\frac{\partial ^{2}\Psi ^{W}}{\partial t^{a}_{\bold {H}}\partial t_{\bold {H}}^{b}}\in \hbar ^{-2}L([\widehat{\gamma }])\]
 Therefore there exist \( \left( A^{(-1)}\right) _{ab}^{c}(t_{\bold {H}}),\, \left( A^{(0)}\right) _{ab}^{c}(t_{\bold {H}})\in \mathbb {C}[[t_{\bold {H}}]] \)
such that 
\[
\frac{\partial ^{2}\Psi ^{W}}{\partial t^{a}_{\bold {H}}\partial t_{\bold {H}}^{b}}-\hbar ^{-1}\sum _{c}\left( A^{(-1)}\right) _{ab}^{c}(t_{\bold {H}})\frac{\partial \Psi ^{W}}{\partial t_{\bold {H}}^{c}}-\sum _{c}\left( A^{(0)}\right) _{ab}^{c}(t_{\bold {H}})\frac{\partial \Psi ^{W}}{\partial t_{\bold {H}}^{c}}\in L([\widehat{\gamma }])\]
 On the other hand, 
\[
\frac{\partial ^{2}\Psi ^{W}}{\partial t^{a}_{\bold {H}}\partial t_{\bold {H}}^{b}},\frac{\partial \Psi ^{W}}{\partial t_{\bold {H}}^{a}}\in L_{W}\widehat{\otimes }\mathbb {C}[[t_{\bold {H}}]]\]
 Therefore 
\[
\frac{\partial ^{2}\Psi ^{W}}{\partial t^{a}_{\bold {H}}\partial t_{\bold {H}}^{b}}=\hbar ^{-1}\sum _{c}\left( A^{(-1)}\right) _{ab}^{c}(t_{\bold {H}})\frac{\partial \Psi ^{W}}{\partial t_{\bold {H}}^{c}}+\sum _{c}\left( A^{(0)}\right) _{ab}^{c}(t_{\bold {H}})\frac{\partial \Psi ^{W}}{\partial t_{\bold {H}}^{c}}\]
 The coordinates \( \{t_{W}\} \) were chosen so that \( \Psi ^{W}\, \text {mod}\, \hbar ^{-1}L_{W}:\mathcal{M}\rightarrow L_{W}/\hbar ^{-1}L_{W} \)
is linear in \( \{t_{W}\} \). Therefore in these coordinates \( \frac{\partial ^{2}\Psi ^{W}}{\partial t^{a}_{W}\partial t_{W}^{b}}=0\, \text {mod}\, \hbar ^{-1}L_{W} \)
and \( \left( A^{(0)}\right) ^{c}_{ab}=0 \). 
\end{proof}
\begin{cor}
One has for one-form \( A=\sum _{a}A_{ab}^{c}(t_{W})dt^{a}_{W} \) : 
\begin{equation}
\label{dA[AA]}
dA=0,\, \, \, \, [A,A]=0
\end{equation}

\end{cor}
The map \( \Psi ^{W} \) induces from the Poincare pairing a bilinear form on
the tangent sheaf of \( \mathcal{M} \) ( \cite{B2} formula (5.59)). If \( W \)
is isotropic with respect to the Poincare pairing: 
\[
(\alpha ,\beta )=0,\, \, \text {for}\, \, \alpha \in W_{\leq r},\beta \in W_{\leq -r+1}\]
 then the form induced on \( T\mathcal{M} \) via \( \Psi ^{W} \)is written
as \( \left\langle \frac{\partial }{\partial t_{W}^{a}},\frac{\partial }{\partial t_{W}^{b}}\right\rangle =\eta _{ab}\hbar ^{n-2} \),
where \( \eta _{ab} \) is graded symmetric and non-degenerate. The series \( A_{ab}^{c}(t_{W}) \)
give structure constants of commutative associative multiplication on the fibers
of \( T\mathcal{M} \) (\( [A,A]=0 \)). They define together with \( \eta _{ab} \)
a quasi-homogeneous solution to WDVV-equation (\cite{B2} theorem 5.6) on (formal)
neighborhood of zero in \( Gr\, W \) identified with \( \mathcal{M} \).

\begin{rem}
Similarly one can construct solutions to WDVV-equation starting from a of kind
of abstract semi-infinite variation of Hodge structure of Calabi-Yau type. More
precisely it is the data \( (\mathcal{G},\Omega [\hbar ]) \) where \( \mathcal{G} \)
is an \( L_{\infty }- \)algebra and \( \Omega [\hbar ] \) is \( \mathbb {C}[\hbar ]- \)linear
\( L_{\infty }- \)module over \( \mathcal{G}\otimes \mathbb {C}^{\hbar }[\xi ] \),
here \( \mathbb {C}^{\hbar }[\xi ] \) is the differential graded commutative
algebra over \( \mathbb {C}[\hbar ] \) with generator \( \xi ,\, \xi ^{2}=0 \),
\( \text {deg}\, \xi =-1 \) and the differential \( d=\hbar \frac{\partial }{\partial \xi } \).
This data define local system \( \mathcal{H}^{\hbar } \) over the moduli space
associated with \( \mathcal{G} \) and the family of semi-infinite subspaces
\( L(t)\subset \mathcal{H}^{\hbar }(t) \) \( t\in \mathcal{M} \), which satisfy
analog of Griffiths transversality condition (\ref{eq:grtransv}). We assume
that moduli space \( \mathcal{M} \) has a smooth subspace \( \widetilde{\mathcal{M}} \)
such that \( L(t)|_{\widetilde{\mathcal{M}}} \) is described by a flat \( \mathcal{O}_{\widetilde{\mathcal{M}}}- \)
module and that Calabi-Yau type condition holds on cohomology level: there exists
\( [\Omega _{0}]\in L(0)/\hbar L(0) \), such that the restriction of symbol
of Gauss-Manin connection (\ref{symbolGMgen}) to \( [\Omega _{0}] \) gives
an isomorphism \( T_{0}\widetilde{\mathcal{M}}\simeq \hbar ^{-1}L(0)/L(0) \).
Then to any subspace \( S\subset \mathcal{H}_{0}^{\hbar } \), \( \hbar ^{-1}S\subset S \)
transversal to \( L(0) \) one can associate a solution \( A^{c}_{ab}(t) \)
to equations (\ref{dA[AA]}). Assume in addition that the module \( \Omega [\hbar ] \)
is equipped with non-degenerate invariant pairing \( G(\hbar a,b)=-G(a,\hbar b)=\hbar G(a,b) \)
which induces pairing on \( L(0) \) with values in \( \hbar ^{D}\mathbb {C}[[\hbar ]] \)
for some \( D\in \mathbb {Z} \). Then for isotropic \( S \) (\( G|_{\hbar S}\in \hbar ^{D}\mathbb {C}[[\hbar ^{-1}]] \))
a flat metrics is induced on \( \widetilde{\mathcal{M}} \) compatible with
multiplication defined by \( A^{c}_{ab}(t) \), i.e. one gets a solution to
WDVV equation. 
\begin{rem}
Using the previous remark one can repeat the above construction of solutions
to WDVV-equation starting from differential graded Lie algebra of Hochcshild
cochains describing \( A_{\infty }- \)deformations of \( A_{\infty }- \)category
\( D^{b}Coh(X) \) and the dg-module describing variations of semi-infinite
Hodge structure associated with local system of periodic cyclic homology. Consequently,
the above formal power series \( A_{ab}^{c}(t_{W}) \) can be constructed entirely
in terms of \( A_{\infty }- \)category \( D^{b}Coh(X) \). 
\end{rem}
Rational Gromov-Witten invariants of Calabi-Yau manifold \( Y \) are encoded
in the series \( C_{ab}^{c}(t_{H},q)^{Y} \) of structure constants of (big)
quantum cohomology (see \cite{KM}). In \cite{B2} we made a conjecture which
can be now reformulated using the fact that the power series \( \Pi ^{W} \)
introduced above coincides with analogous power series from \cite{B2}. The
conjecture states that rational Gromov-Witten invariants of \( Y \) coincide
with Taylor coefficients in the series 
\[
\sum _{\alpha }\frac{\partial ^{2}\Pi ^{W,\alpha }}{\partial t^{a}_{W}\partial t^{b}_{W}}\left( (\partial \Pi ^{W})^{-1}\right) _{\alpha }^{c}=A_{ab}^{c}(t_{W})^{X}\]
 where \( \Pi ^{W,\alpha }(t_{W})=\int _{\Delta _{\alpha }}\Pi ^{W} \) are
the quantum periods considered as functions on the moduli space \( \mathcal{M} \)
associated with the mirror dual family of Calabi-Yau manifolds \( X \). Some
care is needed here: technically we consider family of expansions of \( A_{ab}^{c}(t_{W},q)^{X} \)
with varying base \( [X_{q}]\in \mathcal{M}^{cs} \), \( [X_{q}] \) is close
to point of maximal unipotent monodromy on the boundary of moduli space of complex
structures on \( X \). This is in agreement with the fact that \( C_{ab}^{c}(t_{H},q)^{Y} \)
are formal power series with coefficients in semi-group ring \( Q[B] \) where
\( B \) is the semi-group of algebraic cycles on \( Y \) modulo numerical
equivalences. Interested reader is referred to \cite{B2} \S 6 for details.
Explicitly, the conjecture now states that the identity 
\begin{equation}
\label{mirrident}
C_{ab}^{c}(t_{H},q)^{Y}=A_{ab}^{c}(t_{W},q)^{X}
\end{equation}
 should hold, where \( W \) is the limiting weight filtration on \( H^{*}(X,\mathbb {C}) \)
associated with a maximal unipotency cusp on moduli space of complex structures
on \( X \) and appropriate choice of identification of affine spaces \( H^{*}(Y,\mathbb {C})=Gr\, W \)
is assumed. We checked in \cite{B2} this conjecture for projective complete
intersections: 
\end{rem}
\begin{thm*}
(\cite{B2} theorem 6.2) The identity \( C_{ab}^{c}(t_{H},q)^{Y}=A_{ab}^{c}(t_{W},q)^{X} \)
holds for projective complete intersection Calabi-Yau varieties and their mirrors. 
\end{thm*}
\begin{rem}
One can show using prop. \ref{Pr:PiW}b that in dimensions three this is equivalent
to the standard predictions from \cite{CGOP} proven in \cite{G}. In \cite{G}
some partial higher dimensional generalization is also proven, which deals with
ordinary periods and the subset of rational Gromov-Witten invariants corresponding
to small quantum cohomology ring. 
\end{rem}

\subsection{Quantum periods and mirror symmetry for Fano manifolds.}

A somewhat similar type of integrals appears in a different setting, namely
as periods associated with \( A_{\infty }- \)deformations of a pair \( (X,f) \)
where \( X \) is a quasi-projective variety and \( f:X\rightarrow \mathbb {A}^{1} \)
is a morphism to affine line. Such objects appear as mirror partners to projective
complete intersection Fano manifolds. We consider such a case in \cite{B4}.
In fact in this case the situation simplifies drastically. The Maurer-Cartan
equation becomes empty, the corresponding moduli space can be identified with
the moduli space of deformations of the function \( f \) and the periods of
the semi-infinite Hodge structures are given simply by oscillating integrals.
We prove in loc.cit. an identity analogous to (\ref{mirrident}) which expresses
total collection of rational Gromov-Witten invariants of \( \mathbb {C}\mathbb {P}^{n} \)
in terms of periods of semi-infinite Hodge structure associated with its mirror
partner.

\end{document}